\documentclass[a4paper,12pt,reqno]{amsart}
\usepackage{amsfonts,amsmath,amssymb,amsthm}
\usepackage{graphicx}
\usepackage{fixmath}
\usepackage{hyperref}

\newtheorem{theo}{Theorem}

\DeclareMathOperator{\F}{\mathcal{F}}

\def\phi{\varphi}

\begin{document}

\thispagestyle{plain}

\title{Labeled trees, maps, and an algebraic identity}

\author{Stephan Wagner}
\address{Stephan Wagner\\
        Department of Mathematical Sciences\\
        Stellenbosch University\\
        Private Bag X1\\
        Matieland 7602\\
        South Africa}
\email{swagner@sun.ac.za}


\begin{abstract}
We give a short and direct proof of a remarkable identity that arises in the enumeration of labeled trees with respect to their indegree sequence, where all edges are oriented from the vertex with lower label towards the vertex with higher label. This solves a problem posed by Shin and Zeng in a recent article. We also provide a generalization of this identity that translates to a formula for the number of rooted spanning forests with given indegree sequence.
\end{abstract}

\maketitle

\section{Introduction}
\label{sec:intro}

In a recent paper \cite{shin2011bijective}, Shin and Zeng consider the enumeration of labeled trees with respect to their indegree sequence, where each edge is oriented towards the vertex whose label is higher. This is called the \emph{local orientation} of a labeled tree in \cite{shin2011bijective}, as opposed to the \emph{global orientation} of a rooted labeled tree, in which all edges are oriented towards the root.

In both settings, one can define the indegree sequence $\lambda = 1^{e_1}2^{e_2} \ldots$ (here, $e_i$ is the number of vertices of indegree $i$) of a rooted tree on $n$ labeled vertices, which is always a partition of $n-1$. A classical refinement (see for example \cite[Corollary 5.3.5]{stanley1999enumerative}) of Cayley's formula for the number of labeled trees states that the number of rooted labeled trees with fixed root $r$ and global indegree sequence $\lambda$ is given by
\begin{equation}\label{eq:numtrees}
\frac{(n-1)!^2}{e_0!(0!)^{e_0}e_1!(1!)^{e_1}e_2!(2!)^{e_2}\ldots},
\end{equation}
where $e_0 = n - e_1 - e_2 - \cdots$. Since the relationship between trees and maps will also be important in our present context, it is worthwhile pointing out that~\eqref{eq:numtrees} can be obtained in a very elegant way by considering functional digraphs, see Section 5.3 of \cite{stanley1999enumerative} for details.

Shin and Zeng now gave a bijective proof that~\eqref{eq:numtrees} is also the number of labeled trees with \emph{local} indegree sequence of type $\lambda$, a result that was previously conjectured by Cotterill \cite{cotterill2011geometry} in the context of algebraic geometry and also proved by Du and Yin \cite{du2010counting}.

On the other hand, one has a generating function for labeled trees by their local indegree (see \cite[Eq. (8)]{remmel2002spanning} and also \cite[Theorem 4]{martin2003factorization}): let $\mathcal{T}_n$ be the set of all labeled trees with labels $1,2,\ldots,n$, and let $\operatorname{indeg}_T(i)$ be the local indegree of $i$ in $T$. Then 

$$\sum_{T \in \mathcal{T}_n} \prod_{i=1}^n x_i^{\operatorname{indeg}_T(i)} = P_n(x_1,x_2,\ldots,x_n) = x_n \prod_{j = 2}^{n-1} (jx_j + x_{j+1} + \cdots + x_n).$$

Combining this with~\eqref{eq:numtrees}, one obtains a curious algebraic identity: by the \emph{type} of a monomial $\mathbf{x}^{\mathbf{a}} = x_1^{a_1}x_2^{a_2}\cdots x_n^{a_n}$, we mean the partition $\lambda = 1^{e_1}2^{e_2}\cdots$ of the monomial's total degree in which $e_i$ is the number of exponents $i$ in the monomial. Let now $\lambda = 1^{e_1}2^{e_2}\cdots$ be any fixed partition of $n-1$, and write $e_0 = n - e_1 - e_2 - \cdots$. Then
\begin{equation}\label{eq:identity}
\sum_{\operatorname{type} \mathbf{x}^{\mathbf{a}} = \lambda} [\mathbf{x}^{\mathbf{a}}] P_n(x_1,x_2,\ldots,x_n) = \frac{(n-1)!^2}{e_0!(0!)^{e_0}e_1!(1!)^{e_1}e_2!(2!)^{e_2}\cdots},
\end{equation}
as was pointed out by Shin and Zeng in the aforementioned paper. They asked for a direct proof of this identity, which will be given here. Indeed, we prove a slightly more general identity: 

\begin{theo}\label{thm:main}
f $\lambda = 1^{e_1}2^{e_2}\cdots$ is any partition of $n-k+1$ and $e_0 = n - e_1-e_2-\cdots$, then
$$\sum_{\operatorname{type} \mathbf{x}^{\mathbf{a}} = \lambda} [\mathbf{x}^{\mathbf{a}}] \prod_{j = k}^n (jx_j + x_{j+1} + \cdots + x_n) = \frac{n!(n-k+1)!}{e_0!(0!)^{e_0}e_1!(1!)^{e_1}e_2!(2!)^{e_2}\cdots}.$$
\end{theo}

Note that identity~\eqref{eq:identity} is the special case $k = 2$ of this theorem.

\section{Proof of the main theorem}

We use the standard notation $[n] = \{1,2,\ldots,n\}$. Let us interpret the left hand side in terms of maps from $[k..n] = [n] \setminus [k-1] = \{k,k+1,\ldots,n\}$ to itself, with the additional requirement that $f(r) \geq r$ for all $r$. Let the set of all such maps be $\F_{k,n}$. Multiply the product out to get
$$\prod_{j = k}^n (jx_j + x_{j+1} + \cdots + x_n) = \sum_{f \in \F_{k,n}} C(f) \prod_{j=k}^n x_j^{|f^{-1}(j)|},$$
where $C(f) = \prod_{j \in \operatorname{Fix} f} j$ is the product over all fixed points of $f$. We can eliminate this coefficient by noting that any map $f \in \F_{k,n}$ can be obtained from a map $g: [k..n] \to [n]$ by the rule $f(r) = \Phi g(r) = \max(g(r),r)$, and it is easy to see that each $f$ is obtained from precisely $C(f)$ different maps $g$. Hence the product can also be written as
$$\prod_{j = k}^n (jx_j + x_{j+1} + \cdots + x_n) = \sum_{g: [k..n] \to [n]} \prod_{j=k}^n x_j^{|(\Phi g)^{-1}(j)|}.$$
Any map $h$ whose domain is $[k..n]$ induces a set partition on $[k..n]$ whose blocks are precisely the different preimages. 
Let the type of a set partition be the number partition whose elements are the block sizes. We have to determine the number of maps $g: [k..n] \to [n]$ such that $\Phi g$ induces a set partition of type $\lambda$.

In fact, we can obtain even more: let $\mathcal{A} = (A_1,A_2,\ldots,A_m)$ be any set partition of $[k..n]$ whose type is $\lambda$ (it follows that $m = e_1 + e_2 + \cdots = n - e_0$). Let $a_i$ be the largest element of $A_i$, and assume that $n = a_1 > a_2 > \cdots > a_m$.

Now suppose that $\Phi g$ induces the set partition $\mathcal{A}$. There is only one possibility for the value of $\Phi g(A_1)$, namely $n$, hence $g(r) = n$ for all $r \in A_1 \setminus \{n\}$, but there are $n$ different choices for $g(n) = g(a_1)$. Likewise, there are $n-a_2$ choices for $\Phi g(A_2)$ (namely $a_2,a_2+1,\ldots,n-1$), but if $\Phi g(A_2) = a_2$, one has $a_2$ choices for  $g(a_2)$, so that one has $n-1$ possibilities altogether. Then there are $n-a_3-1$ choices for $\Phi g(A_2)$, but $a_3$ choices for $g(a_3)$ if $\Phi g(A_3) = a_3$ and thus a total of $n-2$ choices, etc.

If follows that the number of possible maps $g$ for which $\Phi g$ induces a fixed set partition $(A_1,A_2,\ldots,A_m)$ is always
$$n(n-1)\cdots (n-m+1) = \frac{n!}{(n-m)!} = \frac{n!}{e_0!}.$$
Furthermore, there are precisely 
$$\frac{(n-k+1)!}{e_1!(1!)^{e_1}e_2!(2!)^{e_2}\cdots}$$
set partitions of $[k..n]$ whose type is $\lambda$, which completes the proof.

\bigskip

Bearing in mind that the number of rooted forests with $r$ given roots on $[n]$ is precisely $r n^{n-r-1}$, it is natural to ask whether the more general identity has an interpretation for forests. This turns out to be the case: let $\mathcal{T}_{r,n}$ be the set of all rooted labeled forests with labels $1,2,\ldots,n$ and roots $1,2,\ldots,r$. We consider the local orientation on these forests. Again by \cite[Eq. (8)]{remmel2002spanning}, one has

$$\sum_{T \in \mathcal{T}_{r,n}} \prod_{i=1}^n x_i^{\operatorname{indeg}_T(i)} = rx_n \prod_{j = r+1}^{n-1} (jx_j + x_{j+1} + \cdots + x_n).$$

Hence we have the following generalization of the result of Du and Yin:

\begin{theo}
The number of rooted labeled forests with labels $1,2,\ldots,n$ and roots $1,2,\ldots,r$ whose type with respect to the local orientation (defined in analogy to the type of a labeled tree) is a given partition $\lambda = 1^{e_1}2^{e_2}\cdots$ of $n-r$ is given by
$$\frac{r(n-1)!(n-r)!}{e_0!(0!)^{e_0}e_1!(1!)^{e_1}e_2!(2!)^{e_2}\cdots},$$
where again $e_0 = n - e_1-e_2-\cdots$.
\end{theo}
Note that an analogous result holds for the global orientation as well (see \cite[Corollary 5.3.5]{stanley1999enumerative}). It is quite probable that the arguments of Du and Yin \cite{du2010counting} and Shin and Zeng \cite{shin2011bijective} for trees carry over to the case of forests as well.

Let us finally point out a reformulation of Theorem~\ref{thm:main} in terms of symmetric functions. Note that the type of a monomial remains unchanged if the variables are permuted. If we consider the symmetrized polynomial
$$\frac{1}{n!} \sum_{\sigma \in S_n} \prod_{j = k}^n (jx_{\sigma(j)} + x_{\sigma(j+1)} + \cdots + x_{\sigma(n)}),$$
then all monomials of the same type $\lambda$ have to have the same coefficient, which is
$$\frac{n!(n-k+1)!}{e_0!(0!)^{e_0}e_1!(1!)^{e_1}e_2!(2!)^{e_2}\cdots} \Bigg/ \frac{n!}{e_0!e_1!e_2!\cdots}. = \frac{(n-k+1)!}{(1!)^{e_1}(2!)^{e_2}\cdots}$$
(dividing by the number of monomials of type $\lambda$). This means, however, that
\begin{equation}\label{eq:symm}
\frac{1}{n!} \sum_{\sigma \in S_n} \prod_{j = k}^n (jx_{\sigma(j)} + x_{\sigma(j+1)} + \cdots + x_{\sigma(n)}) = 
(x_1+x_2+\cdots+x_n)^{n-k+1}.
\end{equation}
It might be interesting to study more general expressions of the form
$$\frac{1}{n!} \sum_{\sigma \in S_n} \prod_{i = 1}^{\ell} \left( \sum_{j=1}^n a_{ij}x_{\sigma(j)} \right)$$
and see under which conditions the sum simplifies as in~\eqref{eq:symm} (and whether there are any combinatorial interpretations, if this is the case).

\bibliographystyle{abbrv}
\bibliography{treesmaps}
 
\end{document}